\documentstyle [12pt] {article}

\thicklines                         
\def\be{\begin{equation}}
\def\ee{\end{equation}}

\begin{document}

\par
\begin{center}
\vspace{0.3in}
\ \\ [0.3in]
{\Large {\bf Nonstandard Poincare and Heisenberg Algebras}}\\
   [.9in]
Preeti Parashar\\ [.2in]
{\it SISSA, Via Beirut 2-4, 34014 Trieste, Italy}\\ [.1in]
{ABSTRACT}\\[.3in]
\end{center}
\begin{quotation}
New deformations of the Poincare group $Fun(P(1+1))$ and its dual 
enveloping algebra $U(p(1+1))$ are obtained as a contraction 
of the $h$-deformed (Jordanian) quantum group $Fun(SL_h(2))$ and its dual.
A nonstandard quantization of the Heisenberg algebra $U(h(1))$ is also 
investigated.
\\ 
[.8in]
\end{quotation}
\vskip1cm
\centerline{Ref. SISSA: 85/96/FM}
\pagebreak
 
Off late, considerable interest has been generated towards the 
nonstandard quantization of Lie groups and algebras, commonly 
known as  $h$ or Jordanian deformation [1-4].
A peculiar feature is that the corresponding universal $R$ matrix is 
triangular ie $R^{-1} = R_{21}$. Hence it is sometimes also called 
triangular deformation. These deformations were further extended to the 
case of supergroups [5].The contraction method is a useful technique to 
study inhomgeneous groups. This was employed first by Celeghini et al 
[11] to obtain quantization of some nonsemisimple groups. Recently 
attempts have been made to apply it to the Jordanian case [6-8] where the
deformation paprameter $h$ has a dimension , like the $\kappa$
deformation [14].

In this letter we propose to obtain a nonstandard quantization of some of 
the simplest inhomogeneous groups - the $(1+1)$ dimensional Poincare group 
, its dual enveloping algebra, and the Heisenberg algebra $U(h(1))$ via a 
contraction of $Fun(SL_h(2))$ and its dual $U_h(sl(2))$. 
The $3$-dim Heisenberg algebra is further extended to $4$ dimensions. 
Another deformation of the 2-dimensional Poincare group can be found in 
[8,9] which was obtained by simultaneously contracting the deformation 
parameter $h$. Our purpose here is to introduce a scaling of the 
generators in such a way that $h$ remains unscaled.\\

$Fun(SL_h(2,R))$ is generated by the matrix T        
\begin{equation}
T=\pmatrix{a&b\cr
c&d\cr}
\end{equation}
modulo the relations 
\begin{equation}\label{comr}
\begin{array}{lll}
{} [c,a] = hc^2 \;, & [b,a] = h - ha^2 \;, & [a,d] = hac - 
hdc\;, \\   
{} [c,d] = hc^2 \;, &[b,d] = h - hd^2 \;, & [c,b] = hac + 
hcd \;.
\end{array}
\end{equation}

The Hopf algebra structure is given by the following
co-product , co-unit and antipode (co-inverse)
\begin{equation}
{\bigtriangleup}(T_{ij}) = T_{ik} \otimes T_{kj} , ~~ i,j = 1,2,
\end{equation}
\begin{equation}
\varepsilon(T_{ij}) = \delta_{ij} ,
\end{equation}

\begin{equation}
S(T)=\pmatrix{d-hc& ~~-b+ha-hd+h^2c\cr 
-c& ~~a+hc\cr}
\end{equation}
and the determinant 
\begin{equation}
D_h = ad - bc -hac = 1
\end{equation}

The universal enveloping algebra $U_h(sl(2))$ intoduced by Ohn [4] is
generated by $J_+$, $J_-$ and $J_3$ with the following Hopf structure
\begin{equation}
[ J_3 , J_+ ] = {2 \sinh(hJ_+)\over h} ,~~
[ J_3 , J_- ] = - \{J_- , \cosh(hJ_+)\} ,~~
[ J_+ , J_- ] = J_3
\end{equation}
\begin{eqnarray}
{\bigtriangleup}(J_+)&=& J_+ \otimes I + I \otimes J_+\nonumber \\
{\bigtriangleup}(J_-)&=&J_- \otimes e^{hJ_+} + e^{-hJ_+} \otimes 
J_-\nonumber \\
{\bigtriangleup}(J_3)&=&J_3 \otimes e^{hJ_+} + e^{-hJ_+} \otimes J_3
\end{eqnarray}                                                      
\begin{equation}
\varepsilon(J_+) = \varepsilon(J_-) = \varepsilon(J_3) = 0
\end{equation}
\begin{equation}
S(J_+) = - J_+ ,~~ S(J_-) = - e^{hJ_+}J_-e^{-hJ_+} ,~~ 
S(J_3) = - e^{hJ_+}J_3e^{-hJ_+}
\end{equation}

The duality between $Fun(SL_h(2))$ and $U_h(sl(2))$ is given by \\
\begin{equation}
\langle J_+ , T \rangle = \pmatrix{0&1\cr
0&0\cr},~~
\langle J_- , T \rangle = \pmatrix{0&0\cr
1&0\cr},~~
\langle J_3 , T \rangle = \pmatrix{1&0\cr
0&{-1}\cr}
\end{equation}
We now apply the following transformation
\begin{equation}
\alpha = a ,~~ \beta = b ,~~ \gamma = {\epsilon}^{-1} c ,~~ \delta = d
\end{equation}

Propostion $1$ : In the limit $\epsilon \rightarrow 0$ these new generators
define an algebra of functions on the nonstandard two dimensional Poincare 
group $P_h(1+1)$. The corresponding algebra and coalgebra is as follows:
\begin{equation}\label{pcom}
\begin{array}{ll}
{} [ \gamma,\alpha ] = 0 \;,& \quad [ \beta,\alpha ] = h - h{\alpha}^2 
\;, \\
{} [ \alpha, \delta ] = 0 \;,& \quad [\beta, \delta ] = h - h{\delta}^2 
\;, \\
{} [ \gamma,\delta ] = 0 \;,& \quad [ \gamma, \beta ] = h\alpha \gamma + 
h\gamma \delta \;,
\end{array}
\end{equation}
\begin{eqnarray}
{\bigtriangleup}(\alpha) &=& \alpha \otimes \alpha,~~~
{\bigtriangleup}(\beta) = \alpha \otimes \beta + \beta \otimes 
\delta,\nonumber \\
{\bigtriangleup}(\gamma) &=& \gamma \otimes \alpha + \delta \otimes 
\gamma,~~~ {\bigtriangleup}(\delta) = \delta \otimes \delta 
\label{gpc}
\end{eqnarray}
\begin{equation}
\varepsilon(\alpha) = \varepsilon(\delta) = 1 , ~~ \varepsilon(\beta) =
\varepsilon(\gamma) = 0
\end{equation}
\begin{equation}
S(\alpha) = \delta ,~~ S(\beta) = -\beta + h\alpha - h\delta ,~~
S(\gamma) = -\gamma ,~~ S(\delta) = \alpha
\end{equation}
and the determinant
\begin{equation}
D_h =  \alpha \delta = 1
\end{equation} 
There are only three independent generators due to the determinant
condition. It is interesting that the expressions for coproduct and counit
of the nonstandard quantum group $Fun(P_h(1+1))$ coincide with that of the
standard Euclidean group $Fun(E_q(2))$ [12]. 

Now by duality, we scale the generators of the enveloping algebra
\begin{equation}
P_+ = J_+ , ~~ P_- = \epsilon J_- ,~~ K = J_3/2
\label{scal}
\end{equation}

Proposition $2$ : $K$, $P_+$ and $P_-$  generate the boosts and the 
translations along the light-cone and define the universal enveloping 
algebra $U_h(p(1+1))$ which is isomorphic to the inhomogeneous algebra
$U_h(iso(1,1))$. \\
This new quantum algebra obeys the following relations which can be deduced
from those of $(U_h(sl(2))$ by the above contraction:
\begin{equation}
[ K , P_+ ] = {\sinh(hP_+)\over h} ,~~
[ K , P_- ] = - P_- \cosh(hP_+) ,~~
[ P_+ , P_- ] = 0
\end{equation}
\begin{eqnarray}
{\bigtriangleup}(P_+)&=& P_+ \otimes I + I \otimes P_+\nonumber \\
{\bigtriangleup}(P_-)&=&J_- \otimes e^{hP_+} + e^{-hP_+} \otimes
P_-\nonumber \\
{\bigtriangleup}(K)&=&K \otimes e^{hP_+} + e^{-hP_+} \otimes K
\label{cop}
\end{eqnarray}  
\begin{equation}
\varepsilon(P_+) = \varepsilon(P_-) = \varepsilon(K) = 0
\end{equation}
\begin{equation}
S(P_+) = - P_+ ,~~ S(P_-) = - P_- ,~~
S(K) = - K + \sinh(hP_+)
\end{equation}

It can be easily checked that the algebras defined in 
Propostions $1$ and $2$ satisfy all the Hopf axioms.
The above Hopf algebra is invariant under the automorphism
$( K, P_+, P_-, h ) \rightarrow ( K, -P_+, -P_-, -h )$.
It was mentioned in [7] that the the nonstandard $(1+1)$ Poincare group
and its dual are isomorphic to each other as Hopf algebras
and one can be obtained from the other by a suitable transformation. It 
should be pointed out that they are isomorphic only at the algebra level. 
At the level of coalgebra  they are different since the 
coproducts are not identical after applying the transformation given in
[7].

The Casimir Operator of $U_h(sl(2))$ is [13 ]
\begin{equation}
C = {{J_3}^2\over 2} + {\sinh(hJ_+)J_-\over h} + {J_- \sinh(hJ_+)\over h}
+ {\cosh^2(hJ_+)\over 2}
\end{equation}
Applying the contraction (\ref{scal}), and, after a suitable renormalization
it leads to the Casimir for $U_h(p(1+1))$
\begin{equation}  
C = {2 P_- \sinh(hP_+)\over h}
\end{equation}
One can see that this is central, and, reduces to the casimir operator of
the classical 2-dim Pioncare algebra as $h \rightarrow 0$.

The universal R-matrix for $U_h(sl(2))$ is given by [9,10,15]
\begin{equation}
R = \exp{{h{\bigtriangleup}J_+\over \sinh(h{\bigtriangleup}J_+)}[J_3 \otimes
\sinh(hJ_+) - \sinh(hJ_+) \otimes J_3]}
\label{unrm}
\end{equation}
In terms of the new generators this becomes
\begin{equation}
R = \exp{{h{\bigtriangleup}P_+\over \sinh(h{\bigtriangleup}P_+)}[K \otimes 
\sinh(hP_+) - \sinh(hP_+) \otimes K]}
\end{equation}
which may be regarded as the universal R-matrix for $U_h(p(1+1))$.

Now we shall obtain a deformation of the 3 dimensional 
Heisenberg algebra $U(h(1))$ again by contracting the algebra of $U_h(sl(2))$
but now with a different scaling. 
Define new generators 
\begin{equation}
A = J_+ ,~~~ A^+ = \epsilon J_- ,~~~ H = \epsilon J_3
\end{equation}
which satisfy the following commutation relations
\begin{equation}
[H, A] = 0 ,~~ [H, A^+] = 0 ,~~ [A, A^+] = H
\end{equation}

Proposition $3$ : $A$ (annihilation operator), $A^+$ (creation operator)
and $H$ (cartan) span a 3-dim deformed heisenberg algebra $U(h_h(1))$.  

Note that the above algebra coincides with its classical counterpart.
However the quantum nature is manifested in the following coalgebra 
structure:
\begin{eqnarray}
{\bigtriangleup}(A)&=& A \otimes I + I \otimes A\nonumber \\
{\bigtriangleup}(A^+)&=&A^+ \otimes e^{hA} + e^{-hA} \otimes
A^+\nonumber \\
{\bigtriangleup}(H)&=&H \otimes e^{hA} + e^{-hA} \otimes H
\end{eqnarray}
\begin{equation}
\varepsilon(A) = \varepsilon(A^+) = \varepsilon(H) = 0
\end{equation}
\begin{equation}
S(A) = - A ,~~ S(A^+) = - e^{hA}A_+e^{-hA} ,~~
S(H) = - H
\end{equation}

In the $q$ case the generators $A$ and $A^+$ were taken to be
hermitian conjugate of each other contrary to the present situation.
The reason for this is that there is no $\ast$ structure defined on the 
Jordanian quantum group and hence does not have a compact form.

We now pass on to investigate the case of 4-dimensional heisenberg algebra.
This is done by central extension i.e.  $sl_h(2) \otimes u(1)$.
Introduce $K$ as the new $u(1)$ generator and apply the transformation
\begin{equation}
A = J_+ , ~~ A^+ = \epsilon J_- ,~~ N = -{J_3\over 2} + {H\over 
2\epsilon} ,~~ H = K
\end{equation}

Proposition 4 : The above transformation in the singular limit
$\epsilon \rightarrow 0$ yields a nonstandard quantization of the 4-dim 
heisenberg algebra which is not a Hopf algebra.

The deformed algebra is
\begin{equation}
[A,A^+] = H , ~~ [N,A] = -{\sinh(hA)\over h} ,~~[N,A^+] = {1\over 2}\{ A^+, 
\cosh(hA) \} ,~~[H, .] = 0
\end{equation}
The coproduct, counit and antipode are given by
\begin{eqnarray}
{\bigtriangleup}(A)&=& A \otimes I + I \otimes A\nonumber \\
{\bigtriangleup}(A^+)&=&A^+ \otimes e^{hA} + e^{-hA} \otimes
A^+\nonumber \\ 
{\bigtriangleup}(H)&=&H \otimes e^{hA} + e^{-hA} \otimes H\nonumber \\
{\bigtriangleup}(N)&=&N \otimes e^{hA} + e^{-hA} \otimes N
\end{eqnarray}
\begin{equation}
\varepsilon(A) = \varepsilon(A^+) = \varepsilon(H) = \varepsilon(N) = 0
\end{equation}
\begin{equation}
S(A) = - A ,~~ S(A^+) = - e^{hA}A_+e^{-hA} ,~~
S(H) = - H , ~~S(N) = - e^{hA}Ne^{-hA}
\end{equation} 

Note that however,the generators $A$ and $N$ form a Hopf subalgebra of 
the above. Infact it turns out that it might not be possible to obtain an 
$h$ Hopf deformation of the above oscillator algebra in the framework of the 
contraction technique.  
 
We remark in passing that it is not straightforward to obtain the universal 
R-matrix for the heisenberg algebra by contracting (\ref{unrm}), as it does 
not involve central factors which could be neglected even if they were 
singular (unlike the case for $q$ deformation [11] ). It would be 
interesting to obtain an expression for such an R-matrix. \\

\vskip.3cm It is a pleasure to thank Dr. B. Jurco for discussions.\\ 
\vskip.3cm
Note added:
After the submission of this work, the paper [16] was brought to my 
attention where the nonstandard Hopf deformation of the $4$-dim oscillator 
algebra was obtained by a different method.

\newpage
{\bf References}
\newcounter{00001}
\begin{list}
{[~\arabic{00001}~]}{\usecounter{00001}
\labelwidth=1cm\labelsep=.5cm}

\item Demidov E.E, Manin Yu.I, Mukhin E.E and and Zhdanovich D.V, {\it Prog.
	Theor. Phys. Suppl.} {\bf 102}(1990)203.
\item Ewen H, Ogievetsky O and Wess J, {\it Lett. Math. Phys.} {\bf 22}(1991)
  297.
\item Zakrzewski S, {\it Lett. Math. Phys.} {\bf 22}(1991)287.
\item Ohn Ch, {\it Lett. Math. Phys.} {\bf 25}(1992)85.
\item Dabrowski L and Parashar P, {\it Lett. Math. Phys.} (in print) 
(q-alg/9511009).
\item Zakrzewski S, { Piosson spacetime symmetry and corresponding 
elementary systems },{\it Quantum Symmetries, Proc. of the II Int. Wigner 
Symposium } (World Scientific, 1993)111.
\item Khorrami M, Shariati A, Abolhassani M.R and Aghamohammadi A, 
{\it Mod. Phys. Lett. A} {\bf A10}(1995)873.
\item Ballesteros A, Herranz F.J, del Olmo M.A, Perena C.M and Santander M,
{\it Jour. Phys. A } {\bf 28}(1995)7113.
\item Shariati A, Aghamohammadi A, Khorrami M , {\it Mod. Phys. Lett.A} 
{\bf 11}(1996)187.
\item Ballesteros A and Herranz F.J, {\it Jour. Phys. A }
{\bf A29}(1996) (to appear) (q-alg/9604008).
\item Celeghini E, Giachetti R, Sorace E and Tarlini M, {\it J. Math. Phys.}
{\bf 32}(1991)1155 ; {\bf 31}(1990)2548.
\item Schupp P, Watts P and Zumino B, {\it Lett. Math. Phys.} {\bf 
25}(1992)139.
\item Ballesteros A, Herranz F.J, del Olmo M.A and Santander M,
{\it Jour. Phys. A} {\bf A28}(1995)941.
\item Lukierski J, Nowicki A and Ruegg H, {\it Phys. Lett.} {\bf 
B293}(1992)344.
\item Vladimirov A.A, {\it Mod. Phys. Lett. A} {\bf 8}(1993)2573. 
\item Ballesteros A, Herranz F.J, {\it Jour. Phys. A} (to appear) 
(q-alg/9602029).
\end{list}
\end{document}